\documentclass[11pt, oneside]{article}   	
\usepackage{amsmath}
\DeclareMathOperator{\maj}{maj}
\DeclareMathOperator{\inv}{inv}

\usepackage{geometry}                		
\usepackage{graphicx}				
\usepackage{amssymb}

\title{A Simple Proof that Major Index and Inversions are Equidistributed}
\author{Michael J. Collins\\Daniel H. Wagner Associates\\\texttt{mjcollins10@gmail.com}}

\begin{document}

\maketitle
\begin{abstract}
We present a short proof of MacMahon's classic result that the number of permutations with $k$ inversions equals the number whose major index (sum of positions at which descents occur) is $k$. 
\end{abstract}

\section{Introduction}
Let $p=p_0p_1\cdots p_{n-1}$ be a permutation of $[n] = \{0,1,\cdots n-1\}$. A \emph{descent} of $p$ is an index $i$ at which $p_{i-1}>p_i$, and an \emph{inversion} of $p$ is a pair of indices $i < j$ with $p_i > p_j$. Define $\inv(p)$ to be the number of inversions in $p$, and define $\maj(p)$, the ``major index" of $p$, to be the sum of all descent positions (so for instance $\maj(241350) = 2+5 = 7$). See (\cite{Bona, Stanley}) for other standard definitions and results regarding permutations.

MacMahon (\cite{MacMahon}) proved that $\inv$ and $\maj$ are \emph{equidistributed}: the number of length-$n$ permutations with $\inv(p)=k$ equals the number of such permutations with $\maj(p)=k$. This common value is denoted $b(n,k)$. MacMahon originally proved this by showing that the generating functions coincide, and Foata (\cite{Foata}) gave a bijective proof; in this note we present a simpler proof.

\section{Proof of Equidistribution}
An \emph{inversion table} of length $n$ is an $n$-tuple of nonnegative integers $(a_0, a_1, \cdots a_{n-1})$ such that $a_j \leq j$ for all $j$. Clearly there are $n!$ inversion tables, and each represents a distinct permutation of $[n]$ as follows: starting with the empty permutation, we repeatedly insert $j$ so that it will have $a_j$ items to its right\footnote{to simplify our presentation we have reversed the usual convention which would have $a_j \leq n-j$}. For instance the inversion table $(0,1,0,3,3)$ yields the permutation $31402$, building it as 
\begin{equation}
\begin{array}{cc}
  &0    \\
  &10    \\
   &102   \\
  &3102  \\
   &34102 
\end{array}
\end{equation}
The insertion of $j$ creates $a_j$ inversions, proving the well-known result that $b(n,k)$ is the number of inversion tables whose elements sum to $k$. To prove equidistribution, we reinterpret $(a_0, a_1, \cdots a_{n-1})$ as meaning repeated insertion of $j$ at a position that will increase the major index by $a_j$. Finding such a position is always possible. For instance (using boldface to emphasize descents) we have $\maj(24\mathbf{1}35\mathbf{0}) = 2+5 = 7$, and the possibilities for insertion of $6$ are:
\begin{equation}
\begin{array}{ccc}
  \maj(24\mathbf{1}35\mathbf{0}6) =& 7+0 =& 2+5   \\
  \maj(24\mathbf{1}356\mathbf{0}) =& 7+1 =& 2+6   \\
  \maj(24\mathbf{1}36\mathbf{5}\mathbf{0}) =& 7+6 =& 2+5+6 \\
  \maj(24\mathbf{1}6\mathbf{3}5\mathbf{0}) =& 7+5 =& 2 + 4 + 6 \\
  \maj(246\mathbf{1}35\mathbf{0}) =& 7+2 =& 3+6 \\
  \maj(26\mathbf{4}\mathbf{1}35\mathbf{0}) =& 7+4 =& 2+3+6 \\
  \maj(6\mathbf{2}4\mathbf{1}35\mathbf{0}) =& 7+3 =& 1+3+6
\end{array}
\end{equation}
In general, say the $\kappa$ inversions of a permutation of $[j]$ occur at positions $d_\kappa < d_{\kappa-1} < \cdots < d_1$. Inserting $j$ at the rightmost position will not change the major index. Insertion at $d_t$ ($1 \leq t \leq \kappa$) will create no new descents, but the descents at $d_t$ through $d_1$ will be shifted to positions $d_t+1, \cdots d_1+1$, so $\maj$ will increase by $t$. Finally, consider inserting $j$ at the $r^{\mbox{th}}$ position (from the left) which is \emph{not} a descent: if there are $r'$ descents to the left of this position, we create a new descent at $r+r'$ and shift $\kappa-r'$ old descents to the right, increasing the major index by $\kappa+r$. Thus the number of permutations with $\maj(p)=k$ is again the number of inversion tables with entries summing to $k$.

\section{Symmetric Joint Distribution}
Our proof of equidistribution is simpler than Foata's, but the machinery of Foata's proof can be used to prove the stronger result that  $\maj$ and $\inv$ have a \emph{symmetric joint distribution} (\cite{Foata, Stanley}): for any pair of integers $k,k'$ the number of $p$ with $\inv(p)=k, \maj(p)=k'$ equals the number with $\inv(p)=k', \maj(p)=k$. We now note that this result can be stated entirely in terms of inversion tables. 

To do this we define another way to interpret an inversion table $(a_0, \cdots a_{n-1})$ as a way to build a permutation, one which makes the relationship between $\inv$ and $\maj$ more direct.
Now $a_j$ will mean ``put $j-a_j$ in the rightmost position, and increment all other elements which are greater than or equal to $j-a_j$". More formally, if $(a_0 \cdots a_{j-1})$ generates the permutation $p_0\cdots p_{j-1}$ then $(a_0 \cdots a_j)$ generates the permutation $p'$ with $p'_j=j-a_j$ and otherwise
\[
p'_k = p_k + [p_k > j-a_j]
\]
Here we make use of the ``Iverson bracket" notation, where $[S]=1$ if the statement $S$ is true, $0$ if it is false.

In fact this just yields the inverse of the permutation generated by reading $(a_0, \cdots a_{n-1})$ as an inversion table. For instance our previous example of $(0,1,0,3,3)$ now yields the permutation $32401$, building it as 
\begin{equation}
\begin{array}{clc}
  & 0 &   \\
  & 10 &   \\
  &  102 &  \\
  &  2130& \\
  &  32401& 
\end{array}
\end{equation}
At step $j$ we create $a_j$ new inversions; the increments do not change any existing inversions, since a pair $r<s$ is either unchanged or becomes $r+1 < s+1$ or $r < s+1$. Furthermore we create a descent at position $j$ if and only if $a_j > a_{j-1}$ (i.e. if position $j$ is an \emph{ascent} of $a$), and similarly the increments do not destroy or create any descents. So the resulting permutation $p$ has
\begin{equation}
\begin{array}{ccc}
\inv(p) &=& \sum a_j \\
\maj(p) &=& \sum_{a_j > a_{j-1}} j
\end{array}
\end{equation}
Therefore, since $\inv$ and $\maj$ are eqidistributed over permutations, the ``sum of elements" and ``sum of ascent positions" are equidistributed over the set of all inversion tables.

\bibliographystyle{plain}
\bibliography{mahonian}

\end{document}